\documentclass{amsart}
\usepackage{amssymb,stmaryrd}
\usepackage{amsfonts}
\usepackage{amstext}
\usepackage{algorithmic}
\usepackage{algorithm}
\usepackage{young}
\usepackage[vcentermath]{youngtab}
\usepackage{graphicx}
\usepackage{epstopdf}
\usepackage[all]{xy}
\usepackage{enumerate}

\usepackage{color}

\numberwithin{equation}{section}
\newtheorem{theorem}{Theorem}[section]
\newtheorem{lemma}[theorem]{Lemma}
\newtheorem{proposition}[theorem]{Proposition}
\newtheorem{corollary}[theorem]{Corollary}

\theoremstyle{definition}
\newtheorem{definition}[theorem]{Definition}
\newtheorem{example}[theorem]{Example}

\theoremstyle{remark}
\newtheorem{remark}[theorem]{\bf{Remark}}

\newcommand{\R}{{\mathbb{R}}}

\newcommand{\C}{{\mathbb{C}}}
\newcommand{\Z}{{\mathbb{Z}}}

\newcommand{\<}{{\langle}}
\renewcommand{\>}{{\rangle}}

\newcommand{\CC}{{\mathcal{C}}}
\newcommand{\CL}{{\mathcal{L}}}

\newcommand{\isom}{{\cong}}

\renewcommand{\ker}{{\rm{ker}}}

\newcommand{\tens}{\otimes}
\newcommand{\id}{\rm id}

\newcommand{\extd}{{\rm d}}
\newcommand{\del}{{\partial}}
\newcommand{\eps}{\epsilon}

\begin{document}

\title{Lie theory and coverings of finite groups} 

\author{S Majid \& K. Rietsch}
\address{Queen Mary University of London\\
School of Mathematical Sciences, Mile End Rd, London E1 4NS, UK}
\address{Kings College University of London\\
Department of Mathematics\\ London, UK}
\email{s.majid@qmul.ac.uk}\email{konstanze.rietsch@kcl.ac.uk}
\thanks{Konstanze Rietsch was supported by Advanced Fellowship grant,
EP/S071395/1, while this research was being conducted}


 \begin{abstract}
We introduce the notion of an `inverse property' (IP)  quandle $\CC$ which we propose as the right notion of `Lie algebra' in the category of sets. To any IP quandle we construct an associated group $G_\CC$. For a class of IP quandles which we call `locally skew' and when $G_\CC$ is finite we show that  the noncommutative de~Rham cohomology $H^1(G_\CC)$ is trivial aside from a single generator $\theta$ that has no classical analogue. If we start with a group $G$ then any subset $\CC\subseteq G\setminus\{e\}$ which is ad-stable and inversion-stable naturally has the structure of an IP quandle. If $\CC$ also generates $G$ then we show that $G_\CC\twoheadrightarrow G$ with central kernel, in analogy with the similar result for the simply-connected covering group of a Lie group. We prove that this `covering map' $G_\CC\twoheadrightarrow G$  is an isomorphism for all finite crystallographic reflection groups $W$ with $\CC$ the set of reflections, and that $\CC$ is locally skew precisely in the simply laced case. This implies that $H^1(W)=k$ when $W$ is simply laced, proving in particular a conjecture for $S_n$ in \cite{Ma:perm}. We obtain similar results for the dihedral groups $D_{6m}$. We also consider $\CC=\Z P^1\cup\Z P^1$ as a locally skew IP-quandle `Lie algebra'  of $SL_2(\Z)$ and show that $G_\CC\cong B_3$, the braid group on 3 strands. The map $B_3\twoheadrightarrow SL_2(\Z)$ which therefore arises naturally as a covering map in our theory, coincides with the restriction of the usual universal covering map $\widetilde {SL_2(\R)}\to SL_2(\R)$ to the inverse image of $SL_2(\Z)$.  \end{abstract}

\maketitle

\section{Introduction}

In the 1950s when they were students Conway and Wraith apparently introduced the notion of a `rack' as capturing the properties of a conjugacy class in a group without having the group itself. The notion came into its own with applications to knot theory \cite{FenRou} and is also of interest in 2-algebra. By a rack we will mean a set $\CC$ with an operation $\CC\times\CC\to \CC$ which we will denote $(a,b)\mapsto {}^ab$ such that
\begin{equation}\label{rack} {}{}^a({}^bc)={}^{({}^ab)}({}^ac),\quad\forall a,b,c\in\CC\end{equation}
and such that ${}^a(\ )$ is bijective for all $a\in\CC$. If ${}^aa=a$ also holds for all $a\in\CC$  then one has a `quandle'. This data is not, however, enough to reconstruct a group in a natural way and in this paper we add a further condition which we call an inverse property (IP) quandle and which in the conjugacy class case corresponds to being stable under group inversion. In Section 2 we show that every IP-quandle leads to an associated group $G_\CC$ with a preferred presentation. 

The construction here is motivated from noncommutative differential geometry applied to finite groups, as in our recent work \cite{LMR}. From this point of view a `Lie algebra' means the dual of the space of invariant differential 1-forms  on the algebra of functions $k(G)$ on the group. The irreducible such objects in the case of bicovariant differential structures are indeed given by conjugacy classes. More generally any ad-stable subset $\CC\subseteq G\setminus\{e\}$ defines a `Lie algebra' and equivalently a bicovariant calculus. We require $\CC$ to generate the group or equivalently the calculus to be connected. For any quasitriangular Hopf algebra, including the finite group case trivially, an inner differential calculus leads to a braided-Lie algebra\cite{GomMa} and in the present case its Lie bracket becomes conjugation. A noncommutative differential calculus is inner if there is a 1-form $\theta$ such that the exterior derivative can be expressed as  $\extd=[\theta,\ ]$, and this is always the case on a finite group or indeed any finite set. In this context it is natural to require $\CC$ to be closed under inversion so as that the differential calculus is symmetric (the associated graph is bidirectional), so $\CC$ is an IP quandle. We show also in Section 2 that in this situation of a group $G$ and $\CC\subseteq G\setminus\{e\}$ we have
\[ G_\CC\twoheadrightarrow G\]
with central kernel. This situation is very similar to that of a Lie group $G$ where its Lie algebra generates an associated connected and simply connected covering group with central kernel\cite{OV}. 

Groups which can be equipped with `Lie algebras' in the sense above such that $G_\CC=G$ are clearly special and we study such covering or `simply connected discrete Lie groups' further. In Section~3 and show that every finite Weyl group  is `simply connected' in this sense with respect to the conjugacy class of reflections, including $S_n$ via its 2-cycles conjugacy class. We also study other quandles for $S_n$ and we cover dihedral groups $D_{2n}$.

Section~4 looks at a first interesting example of a discrete subgroup of a Lie group, $SL_2(\Z)\subset SL_2(\R)$. Here the usual universal covering map $\widetilde{SL_2(\R)}\to SL_2(\R)$ of $SL_2(\R)$ restricts to a central extension $\widetilde{SL_2(\Z)}\to SL_2(\Z)$, where $\widetilde{SL_2(\Z)}$ is the inverse image of $SL_2(\Z)$. This latter map may be thought of as a discrete version of the topological universal covering map, and moreover $\widetilde{SL_2(\Z)}\cong B_3$, the braid group on three strands (see \cite{PRV}). 
We show that this `discrete universal covering' map 
\[ B_3\twoheadrightarrow SL_2(\Z)\]
is an example of our theory where $\CC=\Z P^1\cup\Z P^1\subset SL_2(\Z)$ is a certain natural choice of `Lie algebra' in our sense with $G_\CC\isom B_3$.

In Section~5 we turn to the noncommutative geometry of finite groups obtained from IP quandles, so $G=G_\CC$.  Although it is not clear how to define `simply connected' in noncommutative geometry, an analogue of the abelianization $H^1$ of the fundamental group is understood. We show that this noncommutative de Rham  cohomology $H^1_{d R}$ of $G_\CC$ is essentially trivial when $\CC$ obeys a certain antisymmetry-type condition which  we call `locally skew'. More precisely in this case
\[ H^1_{d R}(G_\CC)=k\theta, \]
the span of the inner generator  $\theta$  which has no classical analogue. Thus the classical part of the cohomology is trivial and the theorem says that our  `covering group' $G_\CC$ in the locally skew case behaves geometrically like the universal cover of a topological group.  We will have shown in Section~4 that this applies to all simply laced Weyl groups, which now proves a conjecture which the first author proposed some time ago\cite{Ma:perm}, that $H^1_{d R}(S_n)=k\theta$ for the 2-cycles class. Our results settle the conjecture for all finite Weyl groups as well as a similar result $H^1_{d R}(D_{6m})=k\theta$ for the dihedral groups $D_{6m}$. We also show that $\Z P^1\cup \Z P^1$ generating $B_3$ is locally skew, so that the noncommutative cohomology here is also trivial when appropriately defined on this infinite group.

\section{IP quandles and its associated group}

Our background is the theory of braided-Lie algebras\cite{Ma:blie}. A braided Lie algebra is a coalgebra $\CL$ in a braided category equipped with a `Lie bracket' $[\ ,\ ]$ subject to certain axioms. Associated to every braided Lie algebra is an induced map $\tilde\Psi:\CL\tens\CL\to \CL\tens\CL$ which obeys\cite{Ma:sol,Wam} the Yang-Baxter or braid relations, and we say that a braided-Lie algebra is regular when this map $\tilde\Psi$ is an isomorphism. In the case of $\CL=k\CC$ in the category of $k$-vector spaces where $\CC$ is a set and $k\CC$ its associated coalgebra, a regular braided-Lie algebra reduces to a rack. In general, associated to a braided-Lie algebra there is a braided enveloping bialgebra $U(\CL)$ defined by relations of $\tilde\Psi$-commutativity. In the rack case the bialgebra becomes the algebra of a monoid with relations of $\tilde\Psi$ commutativity. Explicitly this monoid is
\[ S_\CC=\<\CC\>/ \  {}^ab a=ab,\qquad \forall a,b\in \CC.\]
Here  $\tilde\Psi(a,b)=({}^ab,a)$ is a map $\CC\times\CC\to \CC\times \CC$ and makes $(\CC,\tilde\Psi)$ into a braided set in the sense of \cite{IvaMa}. As a result the monoid $S_\CC$ has reasonable homological properties as an example of an M3-monoid in the notation of \cite{IvaMa}. Note that it is clear that $\CC\subset S_\CC$, since the monoid has quadratic relations and hence is graded by degree, while $\CC$ has degree~1.

\begin{definition} A left IP-quandle is a set $\CC$ with an operation ${}^ab$ obeying (\ref{rack}) and ${}^aa=a$ for all $a$, together with a permutation of $\CC$ denoted $a\mapsto a^{-1}$ such that
\[  (a^{-1})^{-1}=a,\quad {}^a(b^{-1})=({}^ab)^{-1},\quad  {}^{a^{-1}}({}^{a}b)=b,\quad\forall a,b\in \CC.\]
\end{definition}

Clearly this is necessarily a quandle since the third IP-axiom requires each ${}^a(\ )$ to be bijective. It also follows that ${}^{a^{-1}}a=a$ for all $a\in \CC$. The notion is modelled after the notion of an IP-quasigroup and note that we do not require an identity. 

Next, it is easy to see that if $\CC$ is an IP-quandle then $a^{-1}a$ is central in $S_\CC$ for all $a\in \CC$. Indeed, $a^{-1}ab=a^{-1}({}^ab)a={}^{a^{-1}}({}^ab)a^{-1}a=ba^{-1}a$.

\begin{definition} Let $\CC$ be an IP-quandle. We define the group  $G_\CC$ as $S_\CC$ modulo the further relations $a^{-1}a=e$ for all $a\in \CC$. We say that the IP-quandle is {\em embeddable} if  $\CC\subset S_\CC$ descends to an inclusion $\CC\subseteq G_\CC\setminus\{e\}$. 
\end{definition}

Clearly $aba^{-1}={}^abaa^{-1}={}^ab$ in $G_\CC$ for all $a,b\in\CC$. Hence for an embeddable IP-quandle we have $\CC\subseteq G_\CC\setminus\{e\}$, closed under conjugation by $G_\CC$ and clearly closed under inversion.  
 
 \begin{example} Let $\CC=\{a,b,c,a^{-1},b^{-1},c^{-1}\}$ with IP-quandle structure
 \[ {}^aa={}^bc={}^ca=a,\quad {}^ab={}^bb={}^cb=b,\quad {}^ac={}^ba={}^cc=c\]
 and corresponding properties for the inverses (for example ${}^a(c^{-1})=({}^ac)^{-1}=c^{-1}$ etc.) One can verify that the IP quandle axioms hold. However,  viewed in $G_\CC$ we have
 $a=abb^{-1}={}^ab a b^{-1}=bab^{-1}={}^ba=c$ when viewed in $G_\CC$, so this IP-quandle is not embeddable. \end{example}
 
On the other hand, we will normally be interested in IP-quandles constructed from existing groups and in this case the following applies.

\begin{proposition}\label{coverexact} Let $G$ be a group and $\CC$ an ad-stable generating subset $\CC\subseteq G\setminus \{e\}$ closed under group inversion. Then $\CC$ under conjugation and inversion is an embeddable IP-quandle, the canonical inclusion $i:\CC\hookrightarrow G_\CC \setminus \{e\}$ is a map of IP-quandles and covered by the  canonical group homomorphism $\pi:G_\CC\twoheadrightarrow G$ so that $\pi\circ i=\id$. Moreover, its kernel $N$ defined by
\[ 1\to N\to G_\CC{\buildrel \pi\over \to} G\to 1\]
is a central subgroup of $G_\CC$.   \end{proposition}
\proof  This is straightforward. (1) The relations defining $G_\CC$ already hold in $G$, and $\CC$ generates $G$ so the map which associates to $a\in\CC\subset S_\CC$ its corresponding element $a\in G$,  is a surjection and descends to a group homomorphism $\pi:G\to G_\CC$, as stated. The canonical map $i:\CC\to G_\CC$ likewise descending from $\CC\subset S_\CC$ composes with the surjection $G_\CC\twoheadrightarrow G$ to give the inclusion $\CC\hookrightarrow G$,  and hence is itself an injection. The image of $i$ cannot contain $e\in G_\CC$ as this would map under $\pi$ to $e\in G$ and $e\notin\CC$ when viewed in $G$. Moreover,  $i({}^ab)=i(a)i(b)i(a)^{-1}$, $i(a^{-1})=i(a)^{-1}$ by construction of $G_\CC$.  (2) If $n\mapsto e$  under $\pi$ then for all $a\in \CC$ we have $nan^{-1}\mapsto a$. Hence $nan^{-1}={}^na$ and $a$ agree as elements of $\CC$ and hence of  $G_\CC$. Therefore $n\in N=\ker(\pi)$ commutes with all elements of $\CC\subset G_\CC$. Since $G_\CC$ is generated by $\CC$, the subgroup $N$ is central. \endproof

In view of this we henceforth think of $\CC\subseteq G_\CC\setminus\{e\}$ and $\CC\subseteq G\setminus \{e\}$ as sub IP-quandles without writing the map $i$ explicitly.
 
\begin{remark} In the setting of the previous proposition, $G_\CC$ is universal in the sense that if $\phi: H\twoheadrightarrow G$ is a group homomorphism and $j:\CC\hookrightarrow H\setminus\{e\}$ is an inclusion of IP-quandles such that $\phi\circ j=\id$ then there is an induced group homomorphism $\psi$ such that
\[ \begin{array}{ccc} G_\CC & {\buildrel \pi\over \longrightarrow} & G\\ \psi\searrow & &\nearrow\phi \\ &H &\end{array}\]
commutes. 
\end{remark} \proof  Suppose $\phi:H\twoheadrightarrow G$ and  $j:\CC\hookrightarrow H\setminus\{e\}$ obeys $\phi\circ j=\id$ and $j({}^ab)=j(a)j(b)j(a)^{-1}$. We define $\psi(a)=j(a)$ and extend this as a homomorphism of groups. Here $\psi({}^a b a)=j({}^ab)j(a)=j(a)j(b)=\psi(ab)$ and $\psi(a^{-1})=\psi(a)^{-1}$. We have $\phi(\psi(a))=a$ for all $a\in\CC$. Then $\phi(\psi(a_1\cdots a_n))=\pi(a_1\cdots a_n)$ and since $\CC$ generates $G_\CC$ this means that the diagram commutes. \endproof

\begin{remark} Since we can always take $\CC=G\setminus\{e\}$ among other possibilities, Proposition~\ref{coverexact} tells us that up to a central extension every group is `quadratic' in the sense of a set of generators such that the relations are of the form $aa^{-1}=e$ or quadratic, indeed with the quadratic relations obeying the braid or Yang-Baxter identity. \end{remark}

We will be interested in the following `skew-symmetry' condition. We think of IP quandles obeying this condition as being closest in some sense to a classical Lie algebra:
\begin{definition} Let $\CC$ be an IP-quandle. We say that two elements $a,b\in \CC$ are {\em mutually skew} if  ${}^a b=({}^ba)^{-1}$. We say that $\CC$ is {\em skew} if any two distinct elements are mutually skew. Consider the graph on $\CC$ with an edge whenever two distinct elements are mutually skew. We say that $\CC$ is {\em locally skew} if this graph is connected. \end{definition}

Note that from its definition the notion of mutually skew is a property of the unordered pair, i.e. of the subset $\{a,b\}\subseteq\CC$. Also $\{a,a\}$ are mutually skew iff we have $a=a^{-1}$, while $\{a,a^{-1}\}$ are always mutually skew. 

\begin{lemma}\label{skewbraid} Let $\CC$ be an embeddable IP-quandle and $a,b\in\CC$. The following are equivalent
\begin{enumerate}
\item $a,b$ are mutually skew,
\item ${}^{b^{-1}}({}^a(b^{-1}))=a$,
\item ${}^{({}^{b^{-1}}a)}(b^{-1})=a$,
\item the elements $a,b^{-1}$ obey the braid relation $ab^{-1}a=b^{-1}ab^{-1}$  in $G_\CC$,
\item $a^{-1},b^{-1}$ are mutually skew,
\end{enumerate} and are equivalent to the same with $a,b$ interchanged. In this case
\begin{equation}\label{skewIP} {}^{({}^b a)}(a^{-1})=b\end{equation}
\end{lemma}
\proof (2)-(3) are elementary from the IP quandle axioms. For (4) we rearrange as $bab^{-1}=ab^{-1}a^{-1}$ which in view of the relations in $G_\CC$ becomes ${}^ba={}^a{(b^{-1})}$ and the right hand side is $({}^a b)^{-1}$, so this condition is the same as (1). On the other hand (4) can also be rearranged as $a^{-1}b^{-1}a=b^{-1}ab$ which in view of the relations in $G_\CC$ become ${}^{a^{-1}}(b^{-1})={}^{b^{-1}}a$ and the right hand side $({}^{b^{-1}}(a^{-1}))^{-1}$, so this condition is (5). 
Then if $a,b$ and $a^{-1},b^{-1}$ are both mutually skew we have ${}^{({}^ba)}a^{-1}={}^{({}^ab^{-1})}a^{-1}={}^a({}^{b^{-1}}a^{-1})={}^a({}^{a^{-1}}b)=b$. \endproof

Finally, we have been looking at left quandles. There is a parallel theory of right quandles in which we reflect the diagrams at the braided-Lie algebra level so that in our case a right quandle obeys 
\[ (a^b)^c=(a^c)^{(b^c)}\]
with $(\ )^b$ bijective and $a^a=a$ for all $a,b,c\in \CC$. For a right quandle we have  $\tilde\Psi(a,b)=(b,a^b)$ and the relations of $S_\CC$ are $ab=b(a^b)$. By definition the quandle is IP if  there is an involutive bijection denoted as inversion such that  $(a^{-1})^b=(a^b)^{-1}$ and $(a^b)^{b^{-1}}=a$ for all $a,b\in \CC$. 

\begin{lemma} A left IP quandle $\CC$ defines a right IP quandle on the same set with $a^b=({}^{b^{-1}}(a^{-1}))^{-1}$ such that  $(\ )^{-1}$ extends to an isomorphism of the associated groups. \end{lemma}
\proof  Regarding the data for a left quandle as a right one gives a monoid with the opposite product. The corresponding groups are then isomorphic via group inversion and at the quandle level this reduces to the relationship shown.  \endproof

We can similarly convert from right to left in the IP case. Clearly embeddability in one case is likewise equivalent to embeddabilty in the other.

\section{Finite covering groups}

We say that a group $G$ is a covering group if it has a conjugation-stable inversion-stable generating subset $\CC\subseteq G\setminus\{e\}$, regarded as an IP-quandle, such that $G_\CC\twoheadrightarrow G$ is an isomorphism. This means that such a group has a presentation in the form \[ {}^ab a= ab,\quad aa^{-1}=e\quad \forall a,b\in \CC\]
where ${}^ab=aba^{-1}$ is data extracted from the group but only this data is remembered in the form of an IP-quandle. We show in this section that many familiar finite groups are of this form.

\begin{theorem}  Every finite crystallographic reflection group (i.e. Weyl group)  is a covering group with respect to $\CC$ the set of reflections. Moreover, the IP-quandle of reflections is locally skew {\em iff}  the associated Dynkin diagram is simply laced. \end{theorem}
\proof We choose a root system of Dynkin type appropriate to the Weyl group, i.e. a system of roots $\{\alpha\}$ in $\R^n$ for some $n$ with its standard inner product\cite{Hum}. We also choose simple roots $\{\alpha_i\}$. We do not use the full data here and different data can define the same group. Each root   defines a hyperplane reflection which we write as $r_\alpha(v)=v-  \<\check\alpha,v\> \alpha$ for $ v \in \R^n$, where $\check\alpha$ is the coroot of $\alpha$. The quandle $\CC$ can then be identified with the set of roots modulo sign (or all the positive roots) and the IP-quandle structure and relations of $G_\CC$ are represented by
\[{}^{r_\alpha}r_\beta=r_{r_\alpha(\beta)},\quad  r_{r_\alpha(\beta)}r_\alpha=r_\alpha r_\beta,\quad  r_\alpha^{-1}=r_\alpha.\]
 We use the notation $r_i$ for the reflection associated to the simple root $\alpha_i$. Using the relations of $G_\CC$ only, we verify that we can recover all the Coxeter relations as follows. (1) If $\<\check\alpha_i,\alpha_j\>=0$ then ${}^{r_i}r_j=r_j$ and hence $r_jr_i=r_ir_j$. (2) If $\<\check\alpha_i,\alpha_j\>=\<\check\alpha_j,\alpha_i\>=-1$ then ${}^{r_i}r_j=r_{\alpha_i+\alpha_j}$ and hence $r_{\alpha_i+\alpha_j}r_i=r_ir_j$. By symmetry we have similarly $r_{\alpha_i+\alpha_j}r_j=r_jr_i$. Hence $r_ir_jr_i=r_{\alpha_i+\alpha_j}=r_jr_ir_j$ or $(r_ir_j)^3=1$. (3) If $\<\check\alpha_i,\alpha_j\>=-2$ and $\<\check\alpha_j,\alpha_i\>=-1$ we have ${}^{r_i}r_j=r_{\alpha_j+2\alpha_i}$ while ${}^{r_j}r_{\alpha_j+2\alpha_i}=r_{\alpha_j+2\alpha_i}$, hence $r_{\alpha_j+2\alpha_i}r_i=r_ir_j$ and $r_jr_{\alpha_j+2\alpha_i}=r_{\alpha_j+2\alpha_i}r_j$. Then $r_ir_jr_ir_j=r_{\alpha_j+2\alpha_i}r_j=r_jr_{\alpha_j+2\alpha_i}=r_jr_ir_jr_i$ or $(r_ir_j)^4=1$. (4) If $\<\check\alpha_i,\alpha_j\>=-3$ and $\<\check\alpha_j,\alpha_i\>=-1$ we have ${}^{r_i}r_j=r_{\alpha_j+3\alpha_i}$,  ${}^{r_j}r_{\alpha_j+3\alpha_i}=r_{2\alpha_j+3\alpha_i}$ and ${}^{r_i}r_{2\alpha_j+3\alpha_i}=r_{2\alpha_j+3\alpha_i}$, hence $r_ir_j=r_{\alpha_j+3\alpha_i}r_i$, $r_jr_{\alpha_j+3\alpha_i}=r_{2\alpha_j+3\alpha_i}r_j$ and  $r_ir_{2\alpha_j+3\alpha_i}=r_{2\alpha_j+3\alpha_i}r_i$.  Hence $r_ir_jr_ir_jr_ir_j=r_{\alpha_j+3\alpha_i}r_jr_ir_j=r_jr_{2\alpha_j+3\alpha_i}r_ir_j=r_jr_ir_{2\alpha_j+3\alpha_i}r_j=r_jr_ir_jr_{\alpha_j+3\alpha_i}=r_jr_ir_jr_ir_ir_{\alpha_j+3\alpha_i}=r_jr_ir_jr_ir_jr_i$ or $(r_ir_j)^6=1$. Hence $G_\CC$ recovers the Weyl group. 
 
 In the simply laced case, the graph $\Gamma$ of mutually skew edges contains the Dynkin diagram $D$ itself -- two adjacent nodes are labeled by reflections obeying the braid relations hence by Lemma~\ref{skewbraid} are mutually skew. We show that any reflection is connected by a path to this Dynkin diagram. Indeed, any root can be obtained by the action of the Weyl group from a simple root and hence is of the form $r_{i_k}r_{i_{k-1}}\cdots r_{i_1}(\alpha_i)$, say. We apply these simple reflections in turn and show that the accumulated subgraph  of $\Gamma$ remains connected. Thus, applying $r_{i_1}$ to $D$ gives a translate $r_{i_1}(D)$ of $D$ as a subgraph of $\Gamma$, because the Weyl group acts as a group automorphism by conjugation and hence two vertices are connected after the action of $r_{i_1}$ iff they  were before. The translate is connected to $D$ since at least one node, $\alpha_{i_1}$ itself, is invariant so lies in both $D$ and its translate. We write $\Gamma_1=D\cup r_{i_1}(D)$ for the combined connected subgraph. Now we apply $r_{i_2}$ to $\Gamma_1$ to obtain a subgraph $r_{i_2}(\Gamma_1)$ as its translate. This translate remains connected to $D$ since $\alpha_{i_2}\in D$ is invariant and so also lies in $r_{i_2}(\Gamma_1)$. We write $\Gamma_2=D\cup r_{i_2}(\Gamma_1)$ for the combined connected graph. Iterating this process with $\Gamma_{k+1}=D\cup r_{i_{k+1}}(\Gamma_{k})$ we arrive at a connected subgraph of $\Gamma$ which contains the element of the Weyl group that we began with as well as $D$. Hence the group is locally skew. 
 
Next we claim that if $r_\alpha,r_\beta$ are mutually skew then $\alpha,\beta$ have the same length. Hence in the non-simply laced case the graph of mutually skew edges must have different components and the quandle is not locally skew. Indeed,  $r_{\alpha}$ and $r_{\beta}$ being mutually skew means 
$r_\alpha(\beta)=\pm r_\beta(\alpha)$, i.e. $\beta-\alpha\<\check\alpha,\beta\>=\alpha-\beta\<\check\beta,\alpha\>$ or $\beta-\alpha\<\check\alpha,\beta\>=-\alpha+\beta\<\check\beta,\alpha\>$. If $\alpha\ne\pm\beta$ we conclude in the first case that $r_\alpha(\beta)=\alpha+\beta=r_\beta(\alpha)$ and in the second case that $r_\alpha(\beta)=\beta-\alpha=-r_\beta(\alpha)$. In the first case, since $\alpha=r_\alpha(\beta)-\beta$, we deduce that $r_\beta(\alpha)=r_\beta(r_\alpha(\beta)-\beta)=r_\beta r_\alpha(\beta)+\beta=\alpha+\beta$ hence $r_\beta r_\alpha(\beta)=\alpha$. Hence $\alpha$ has the same length as $\beta$ as reflections do not change length. There is a similar argument in the other case.  \endproof

\begin{example} Let $\CC\subset S_n$ be the conjugacy class of 2-cycles. Clearly 
\[ {}^{(ij)}(kl)=(kl),\quad  {}^{(ij)}(ij)=(ij),\quad  {}^{(ij)}(jk)=(ik)\]
for disjoint $i,j,k$ and $(ij)^{-1}=(ij)$. Hence $\CC$ is a skew IP-quandle if $n\le 3$ and clearly a locally skew one if $n>3$. The relations of $G_\CC$ are 
\[ (ij)(kl)=(kl)(ij),\quad (ik)(ij)=(ij)(jk),\quad (ij)^2=e\]
for disjoint $i,j,k$. These imply the usual Coxeter relations since $(i,i+1)(i+1,i+2)(i,i+1)=(i,i+2)(i,i+1)^2=(i,i+2)(i+1,i+2)^2=(i+1,i+2)(i,i+1)(i+1,i+2)$ using only the relations above. Hence the relations stated form a presentation of $S_n$  and hence $G_\CC=S_n$ as per the theorem.
\end{example}

We conclude with some other examples. The first illustrates that a group can be a covering group with respect to one generating $\CC$ and not with respect to another. We will always write $\Z_m$ for the cyclic group $(\Z/m\Z, +)$. 

\begin{example} Let $\CC\subset S_4$ be the conjugacy class of 4-cycles  as an IP-quandle. It consists of the following elements $e_i$ and their inverses,
\[ e_1=(1234),\quad e_2=(2134),\quad e_3=(1324)\]
and the quandle structure is
\[  {}^{e_i}e_j=e_k^{\eps_{ijk}},\quad  {}^{e_i}e^{-1}_j= e_k^{-\eps_{ijk}},\quad  \forall i\ne j,\]
where $k$ is the other one of the triple $i,j,k$ and $\eps$ is totally antisymmetric with $\eps_{123}=1$. Hence the IP-quandle here is skew. The relations of $G_\CC$ are those expressing inverse and 
\[ e_ie_j=e_k^{\eps_{ijk}} e_i,\quad\forall i\ne j\]
with the conventions above. Within $S_4$ we would have  $(12)=e_3 e_1^2, (23)=e_2e_3^2, (34)=e_3 e_2^2$ and taking these now as definitions we seek to show the relations of $S_4$. The first step is to show that $e_1^2e_2^2=e_3^2$, $e_2^2e_1^2=e_3^{-2}$ and cyclic rotations of these. The proof of the first of these is $e_1^2=e_1e_2e_2^{-1}e_1=e_3e_1e_3e_2^{-1}=e_3^2e_2^{-2}$ while the second is $e_2^2e_1^2=e_2e_1e_3^{-1}e_1=e_1e_3^{-1}e_3^{-1}e_1=e_1e_3^{-1}e_2^{-1}e_3^{-1}=e_3^{-2}$. From these equations we deduce that $e_1^2e_2^4e_1^2=e$ or $e_2^4=e_1^{-4}$ which similarly equals to $e_3^4=e_2^{-4}$. We conclude that $e_i^8=e$. Also using the relations,  the element $e_1^4=e_2^4=e_3^4$ is central so this element  $\zeta$, say, generates a central $\Z_2$ subgroup while the $e_i^2$  generate a copy of the quaternion group $Q_8$. We also have  $e_i^2e_j=e_j^{-1}e_i^2$ and $e_i^2e^{-1}_j=e_je_i^2$ for all $i\ne j$. It is easy to verify that if we quotient by the above $\Z_2$ then the relations $(12)^2=e$ etc and the rest of the Coxeter relations hold. Hence 
\[ 1\to \Z_2\to G_\CC\twoheadrightarrow S_4\]
so that $G_\CC$ is a double cover. Moreover, from $(12)^2=(23)^2=(34)^2=\zeta$ and other relations in $G_\CC$, it is clear that $G_\CC\isom 2S_4^-$, one of the two Schur covers of $S_4$.   \end{example}

This example restricts on the inverse image of the Klein four group $D_4=\Z_2\times\Z_2$ to the quaternion group  $Q_8$ as its double cover. Note that the $n$-cycles form a generating IP quandle in $S_n$ for all $n$ even, and in $A_n$ for all odd $n$. However, we do not know in general the kernels of the maps $G_\CC\twoheadrightarrow S_n$ for even $n$ and $G_\CC\twoheadrightarrow A_n$ for odd $n$ nor if the IP quandles are in general locally skew.

\begin{proposition} The dihedral groups $D_{2n}=\Z_n\rtimes\Z_2$ are  covering groups with respect to their generating sets $\CC=\Z_nx$, where $x$ is the generator of $\Z_2$. Moreover, $\CC$ is skew iff $n=1,3$ and for $n>1$ is locally skew iff $n$ is divisible by 3. \end{proposition}
\proof Here $D_{2n}$ has generators $a,x$ and relations $a^n=x^2=e$, $xa=a^{-1}x$ and $\CC=\Z_nx$ generates (this is a $n$ conjugacy class when $n$ is odd and a sum of two classes when $n$ is even). Writing its members $e_i=a^ix$ where $i=0,\cdots,n-1$ we have 
\[ {}^{e_i}e_j=e_{2i-j},\quad i,j\in \Z_n\]
where the indices are treated modulo $n$. Hence the IP-quandle here is skew iff $n=1,3$. More generally $\{e_i,e_j\}$ are mutually skew iff $n$ divides $3(i-j)$. This means that unless $n=3m$ (so $D_{6m}$) we have no distinct mutually-skew pairs, while for $D_{6m}$ we have for example $\{e_i,e_{i+1}\}$ are mutually skew for all $i\in\Z_n$ and hence $D_{6m}$ is locally skew. The relations of $G_\CC$, where $e_i$ are regarded abstractly as generators, are
\[ e_{i-j}e_i=e_ie_{i+j},\quad e_i^2=e,\quad \forall\ i,j\in\Z_n\]
Conversely given this group $G_\CC$ we let $a=e_1e_0=e_1e_2=\cdots=e_{n-1}e_{n-2}=e_0e_{n-1}$ by repeatedly using the relations. We then have $a^i=e_{i}e_0=e_{i+1}e_1=\cdots =e_{i+n-1}e_{n-1}$ as we prove by induction: suppose true for $a^{i-1}$ then $a^i=aa^{i-1}=(e_{j+1}e_j)(e_j e_{j-i+1})=e_{j+1}e_{j+1-i}$ for all $j$ under our induction hypothesis.  In particular, $a^n=e_0e_0=e$. We also let $x=e_0$ so that $x^2=e$ and $xa=e_0e_1e_0=e_1e_2e_0=a^{-1}x$. Hence all the relations of $D_{2n}$ are recovered and $G_\CC=D_{2n}$. 
\endproof

This includes $D_2=\Z_2$ and $D_4=\Z_2\times\Z_2$ as covering groups with respect to the stated $\CC$. On the other hand it is easy to see that if $G=\Z_n$, $n>2$ and $\CC=\{a,a^{-1}\}$ where $a$ is a generator then 
\[ G_\CC=\Z\twoheadrightarrow\Z_n\]
 as covering group. The same applies for any finite abelian group with factors $\Z_2$ or $\Z_n$, $n>2$ and $\CC$ given by one generator for each cyclic group, giving respectively $\Z_2$ or $\Z$ factors in the covering group. So among finite Abelian groups only $\Z_2^m$, $m\ge 1$ are covering groups. 


\section{Covering group of  $SL_2(\Z)$}

We are not limited to finite groups and quandles and in this section we describe an infinitely generated example based on  $SL_2(\Z)$. We let $\CC=\Z P^1\cup\Z P^1=\{e_{\vec a}, e^{-1}_{\vec a}\}$ where $\vec a\in \{(a,c)\in \Z^2\ |\ (a,c)=1\}/\Z_2$ (coprime integers modulo $\Z_2$ acting by -1 and noting that $(0,1)=(1,0)=1$ so these are included) and
\[ e_{ (a,c)}=\left(
\begin{array}{cc}
 1-a c & a^2 \\
 -c^2 & 1+a c
\end{array}
\right),\quad e_{(a,c)}^{-1}=\left(
\begin{array}{cc}
 1+a c & -a^2 \\
 c^2 & 1-a c
\end{array}
\right).\]
This the union of the conjugacy classes of $e_{(1,0)}=\begin{pmatrix}  1&1\cr 0 & 1\end{pmatrix}$ and $e^{-1}_{(0,1)}=\begin{pmatrix}1&0\cr 1&1\end{pmatrix}$. We have an IP-quandle with 
\[ {}^{e_{\vec a}}e_{\vec b}=e_{\vec a|{{\vec a}\atop{\vec b}}|+\vec b},\quad  {}^{e_{\vec a}}e^{-1}_{\vec b}=e^{-1}_{\vec a|{{\vec a}\atop{\vec b}}|+\vec b},\quad {}^{e^{-1}_{\vec a}}e_{\vec b}=e_{-\vec a|{{\vec a}\atop{\vec b}}|+\vec b},\quad  {}^{e^{-1}_{\vec a}}e^{-1}_{\vec b}=e^{-1}_{-\vec a|{{\vec a}\atop{\vec b}}|+\vec b}\]
and the associated group $G_\CC\twoheadrightarrow SL_2(\Z)$ with relations 
\[ e_{\vec a}e_{\vec b}=e_{\vec a|{{\vec a}\atop{\vec b}}|+\vec b}e_{\vec a}. \]
Here $|\ |$ denotes determinant. It is easy to see that $\{e_{\vec a},e_{\vec b}\}$ are never mutually skew, and $\{e_{\vec a},e^{-1}_{\vec b}\}$ are mutually skew iff either $\vec a=\pm \vec b$ (which means $\{e_{\vec a},e^{-1}_{\vec a}\}$) or $|{\vec a\atop \vec b}|=\pm 1$ (it is always possible to find $\vec b$ for a given $\vec a$ obeying this by the Bezout identity). Each such pair then gives by Lemma~\ref{skewbraid} a homomorphic image of the braid group $B_3$ in $G_\CC$ via $e_{\vec a}$ and $e_{\vec b}$. For example $e_{(1,0)}$, $e^{-1}_{(0,1)}$ are mutually  skew and hence give such a homomorphic image by $e_{(1,0)},e_{(0,1)}$. 

\begin{proposition} The IP quandle $\Z P^1\cup\Z P^1$ is locally skew. Moreover, the associated group $G_\CC$ is isomorphic to $B_3$, the braid group on 3 strands.
\end{proposition}
\proof (1) We construct a path between $e_{\vec a}$ and $e^{-1}_{\vec b}$  or $e_{\vec b}$ by steps of the form $\{e_{\vec a},e^{-1}_{\vec a_1}\}, \{e^{-1}_{\vec a_1},e_{\vec a_2}\},\cdots, \{e_{\vec a_n},e^{-1}_{\vec b}\}$ in the first case and similarly in the other case. We can insert a pair $\{e_{\vec b},e_{\vec b}^{-1}\}$ as necessary to make our sequence of this form. Then by the above observations, we need to find coprime vectors in $\Z^2$ so that adjacent members of  the sequence $\vec a, \vec a_1,\cdots,\vec a_n,\vec b$ have mutual determinant $\pm 1$. This means that from the origin they form a triangle of area 1/2. This can be achieved as follows: form a triangle from the origin with sides $\vec a,\vec b$. We then triangulate the interior of this with vertices taken from $\Z^2$. This can always be achieved by induction (pick any interior vertex to make a triangle and then triangulate the remainder). Then starting at $\vec a$, the desired sequence is provided by the edge vectors in turn rotating about the origin until reaching $\vec b$. Each triangle formed from the adjacent edges has no interior points and hence area 1/2 by Pick's theorem.  (2) Observe first that for any group $G$ and ad-stable inversion stable generating set $\CC\subseteq G\setminus \{e\}$, if a subset $S\subseteq\CC$ generates $G$  as a group and generates $\CC$ as an ad-stable subset then $S$ also generates $G_\CC$. This is because if $a\in \CC$ we realise it as $a=(s_1\cdots s_{n-1})s_n(s_1\cdots s_{n-1})^{-1}= {}^{s_1}({}^{s_2}(\cdots({}^{s_{n-1}}s_n)\cdots))$ for some generators. Viewing this in $G_\CC$ we obtain $a$ as a series of conjugations by the $s_i$. We can then use that $\CC$ generates. Now in our case the elements  $e_{(1,0)}$ and $e^{-1}_{(0,1)}$ generate $SL_2(\Z)$ and hence $G_\CC$ and hence the map from $B_3$ described above is surjective. 
The map is injective because the kernel of the map $B_3\to SL_2(\Z)$ is known to be $\Z$ generated by the element $(b_1b_2)^{6}$ where $b_1,b_2$ are the standard generating braids. Hence any elements  in the kernel of the factoring map $B_3\to G_\CC$ are powers of $(b_1b_2)^{6}$. But  $(e_{(0,1)}e_{(1,0)})^{6}$ has infinite order because among  the $e_{\vec a}$ the relations are all quadratic and hence any product of them can never be the identity.   \endproof

There is an identical story with $PSL_2(\Z)$. Because elements of $\CC$ above all have trace $2$ they are all distinct in $PSL_2(\Z)$, so it has the same `Lie algebra'.

\section{Cohomology $H^1_{d R}(G_\CC)$ }

The general notion of differential graded algebras over algebras is one of the components of noncommutative geometry\cite{Con}. In particular,  left and right translation covariant or `bicovariant' differentials on Hopf algebras have been much studied since their formulation in \cite{Wor}. The specialisation to finite groups is also well-known and we refer to \cite{Ma:cdo,Ma:prim} for an introduction. Briefly, let $G$ be a finite group and $k(G)$ be the algebra of functions. Bicovariant  noncommutative differential calculi $\Omega^1(k(G))$ correspond to Cayley digraphs where arrows are of the form $x\to xa$ for $a\in \CC$,  a fixed ad-stable subset of $G\setminus\{e\}$. Here $\Lambda^1={\rm span}_k\{\omega_a\}_{a\in\CC}$  is the space of left-invariant 1-forms and in terms of the Kronecker delta-functions the basis 1-forms are given by $\omega_a=\sum_{x\in G}\delta_x\extd\delta_{xa}$. We have 
\[ \omega_a.f=R_a(f)\omega_a,\quad \extd
f=\sum_{a\in \CC}\del^a(f)\omega_a,\quad\forall f\in k(G),\]
where the finite partial derivative is $\del^a=R_a-\id$ in terms of right translation $R_a(f)=f((\ )a)$.  The calculus is inner with $\theta=\sum_a\omega_a$ in the sense that $\extd f=\theta f - f\theta$ for all $f\in k(G)$, and is connected iff $ \CC$ is a generating set, which we suppose. It is also natural to suppose $\CC$ is stable under inversion so that the Cayley graph is bidirected (every arrow has an inverse) and we suppose this, as in earlier sections. The stated $\Lambda^1$ extends in a natural way to a braided Hopf algebra $\Lambda$ as the invariant part of an exterior algebra $\Omega=\oplus_i\Omega^i$ where $\Omega^0=k(G)$ and this forms a complex with $\extd^2=0$, the noncommutative de Rham complex. Its cohomology is denoted $H^1_{d R}(G)$ and depends on $\CC$. The relations of the exterior algebra in degree 2 are quadratic and take the form
\[ \sum v_{a,b}\omega_a\wedge \omega_b=0,\quad{\rm iff}\quad \sum v_{a,b}\omega_a\tens \omega_b\in \ker(\id-\tilde\Psi)\]
where $\tilde\Psi(\omega_a\tens \omega_b)=\omega_{aba^{-1}}\tens \omega_a$ is adjoint to the $\tilde\Psi$ on $\CC$ as a right-handed rack. Here $\Lambda^1$ and $k\CC$ are mutually dual.  One has $\extd\omega=\theta\omega+\omega\theta$ for all $\omega\in \Omega^1$ and similarly with graded commutator in all degrees. The higher degrees are similarly determined by $\tilde\Psi$ in a canonical way\cite{Ma:cdo,Wor}.

\begin{theorem}\label{cohtriv} Let $\CC$ be a right IP-quandle and $k$ not have characteristic 2. Then $H^1_{d R}(G_\CC)\supseteq k\theta$, the span of the element $\theta$. Equality holds if $\CC$ is locally skew.
\end{theorem}
\proof Let $\CC$ be a right IP-quandle. We note first that
\[ \tilde\Psi(\omega_a\tens\theta)=\theta\tens\omega_a,\quad \tilde\Psi(\theta\tens\omega_a)=\sum_{b}\omega_{bab^{-1}}\tens\omega_b\]
Hence if $\omega=\sum_a c_a\omega_a$ where $c_a\in k(G)$ and $\extd\omega=0$ in $\Omega^2$, i.e. using the Leibniz rule, 
\[ \sum_{a,b}(\del^b c_a)\omega_b\tens\omega_a+\sum_a c_a(\theta\tens\omega_a+\omega_a\tens\theta)\in \ker(\id-\tilde\Psi).\]
which is equivalent to 
\[ \sum_{b,a}R_bc_a\omega_{bab^{-1}}\tens\omega_b- \sum_a c_a\omega_a\tens\theta=\sum_{a,b}R_b c_a\omega_b\tens \omega_a - \sum_a c_a\theta\tens\omega_a\]
or 
\[ R_b c_{b^{-1}ab}+c_b= R_a c_b + c_a,\quad \forall a,b\in\CC.\]
Now consider 
\[ f(x)= c_{a_1}(e)+c_{a_2}(a_1)+\cdots +c_{a_n}(a_1\cdots a_{n-1}),\quad x=a_1a_2\cdots a_n\]
where we write $x$ as obtained along a path from $e$ in the Cayley graph by repeated multiplication from the right by generators $a_1,a_2,\cdots$. We add up the values of $c$ at each interim vertex in the direction of the next vertex. That this is independent of the path requires invariance under application of the relations in the group $G_\CC$ which are either  the quadratic ones or  express inverses. This means in diagram form:
\[ \includegraphics[scale=.7]{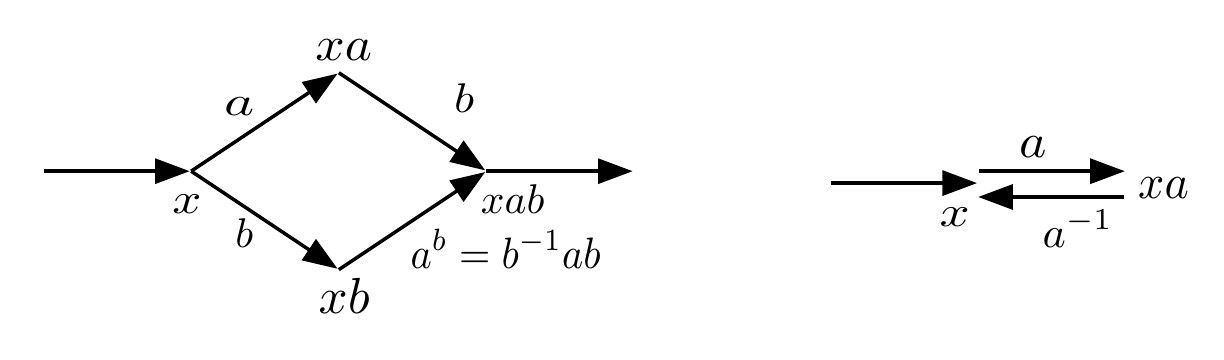}\]
at all $x\in G$ and $a,b\in\CC$. The first requires 
\begin{equation}\label{cocy}c_a(x)+c_b(xa)=c_b(x)+c_{a^b}(xb)\end{equation}
 which is the cocycle condition we found on $c$, while the second 
requires $c(x)_a+c(xa)_{a^{-1}}=0$. If these conditions hold then $f$ by construction obeys $\del^a f=c_a$ and hence $\omega=\extd f$. Conversely, if $\omega$ is exact then these conditions hold. We see in particular that $\theta$ is closed but not exact. More generally, if 
\[ c_a(x)+c_{a^{-1}}(xa)\]
is a constant $\lambda$ independent of $x,a$, then clearly $\bar c_a=c_a-{\lambda\over 2}$ obeys the conditions that we required for $c$ above, i.e. $\omega-{\lambda\over 2}\theta$ is exact. The following lemma then completes the proof. \endproof

\begin{lemma} Let $\CC$ be a right IP-quandle, $\omega=\sum_ac_a\omega_a$  a closed 1-form and $d_a(x)=c_a(x)+c_{a^{-1}}(xa)$. \begin{enumerate}
\item  $d_a(x)=d_{a^b}(xb)$ for all $a,b\in\CC$ and $x\in G_\CC$.
\item If $\{a,b\}$ are mutually skew then $d_a(x)=d_b(x)$ for all $x\in G_\CC$.
\end{enumerate}
\end{lemma}
\proof From (\ref{cocy}) at we also have 
\[ c_{a^{-1}}(xa)+c_b(x)=c_b(xa)+c_{(a^{-1})^b}(xab)\]
and adding this to to (\ref{cocy}) we have an equation which we use to compute
\[ d_{a^b}(xb)=c_{a^b}(xb)+c_{(a^b)^{-1}}(xba^b)=c_{a^b}(xb)+c_{(a^{-1})^b}(xab)=c_a(x)+c_{a^{-1}}(x)=d_a(x)\]
using also the definition of an IP quandle and the relations of the group. Hence $d_a(x)=d_{a^{x^{-1}}}(e)$ where we extend the action of $\CC$ to an action of $G_\CC$. In view of this we will focus on $x=e$ for simplicity. We now return to equation (\ref{cocy}) with $x,a,b$ assigned different values and overall signs as follows:
\begin{eqnarray*} (x,a,b)=(e,u,v):& c_u(e)+c_v(u)&=c_v(e)+c_{u^v}(v)\\
 (x,a,b)=(e,u,(u^{-1})^{v^{-1}}):& c_u(e)+c_{(u^{-1})^{v^{-1}}}(u)&=c_{(u^{-1})^{v^{-1}}}(e)+c_{u^{((u^{-1})^{v^{-1}})}}((u^{-1})^{v^{-1}})\\
 (x,a,b)=(u,v,u^{-1}):& -c_v(u)-c_{u^{-1}}(uv)&=-c_{u^{-1}}(u)-c_{v^{u^{-1}}}(e)\\
 (x,a,b)=(v,u^v,v^{-1}):& c_{u^v}(v)+c_{v^{-1}}(uv)&=c_{v^{-1}}(v)+c_u(e)\\
 (x,a,b)=(uv,u^{-1},v^{-1}):& c_{u^{-1}}(uv)+c_{v^{-1}}(v^{u^{-1}})&=c_{v^{-1}}(uv)+c_{(u^{-1})^{v^{-1}}}(u)
 \end{eqnarray*}
 where we used definition of an IP-quandle and the relations of the group to simplify. We now add all these equations together and cancel to obtain
 \[ d_u(e)=d_v(e)+c_{(u^{-1})^{v^{-1}}}(e)+c_{u^{((u^{-1})^{v^{-1}})}}((u^{-1})^{v^{-1}})-c_{v^{u^{-1}}}(e)-c_{v^{-1}}(v^{u^{-1}})\]
 for all $u,v\in\CC$. If $u^{-1},v^{-1}$ are mutually skew then $(u^{-1})^{v^{-1}}=v^{u^{-1}}$ and if also $u,v$ are mutually skew then  $u^{((u^{-1})^{v^{-1}})}=v^{-1}$ by the right handed version of (\ref{skewIP}). According to Lemma~\ref{skewbraid} this is not a further requirement. Hence in this case $d_u(e)=d_v(e)$.  \endproof

This then completes the proof of the theorem: when $\CC$ is locally skew we conclude that $d_a(e)$ is independent of $a\in\CC$ by part (2) and hence $d_a(x)$ is independent of $a,x$ by part (1).   From Section~3 we see among others that $H^1_{d R}(S_n)=k\theta$ for all $n$ where $S_n$ has its 2-cycles differential calculus and $H^1_{d R}(D_{6m})=k\theta$ with $\CC=\Z_{3m}x$.

If $\CC$ is not locally skew the proof of the theorem and the lemma still tells us that a closed 1-form $\omega=\sum_ac_a\omega_a$ is exact iff $d_a(e)=0$ for all $a$, although not all of these conditions may be independent of each other or of the cocycle condition already imposed. This implies that 
\[ \dim H^1_{d R}(G_\CC)\le |\CC|\]
for any finite IP-quandle and its associated group.  

\begin{corollary} If $\CC$ is a trivial IP-quandle and $k$ does not have characteristic 2 then $H^1_{d R}(G_\CC)=k\CC$ spanned by the basis $\{\omega_a\}$ of left-invariant 1-forms. \end{corollary}
\proof In this case the cocycle condition (\ref{cocy}) is solved by constant coefficients $c_a(x)$ independent of $x$. This is the space of left-invariant 1-forms. On this space the additional condition to be exact is the $|\CC|$ equations $d_a=2c_a=0$. Assuming that the characteristic is not 2 this means that $H^1_{d R}$ contains each $\omega_a$ as representatives. These are already linearly independent hence the bound is saturated.  \endproof

Thus for $D_{4}=\Z_2\times\Z_2$ we have $H^1_{d R}(D_4)=k^2$ spanned by $\omega_0,\omega_1$ associated to the elements $e_0,e_1$ of $\CC$. These are the left-invariant 1-forms and it is known that they provide the cohomology in this and similarly trivial examples. 

\begin{remark} we have discussed only $H^1_{d R}$ above but we remind the reader of the conjecture in \cite{Ma:perm} that the volume dimensions (top degrees) of $\Omega(\C(S_n))$ with $ \CC=\{2-cycles\}$ are the same as the number of indecomposable modules of the `preprojective algebra' of type $sl_n$ (related to the Lusztig canonical basis of $U_q(sl_n)$). One could also expect at the top degree that $H^{top}_{d R}=\C=H^0_{d R}$ by Poincar\'e duality and then the conjecture becomes that this top cohomology dimension is the number of indecomposables. Similarly for the Weyl groups of other Lie algebras in relation to their canonical bases. The work \cite{Ma:perm} showed how a symmetric algebra partner to the left invariant exterior algebra naturally recovers the Fomin-Kirillov algebra that arises in the classical cohomology of flag varieties in type A,  an approach which was extended to general type in \cite{Baz}. The work \cite{Ma:perm} also considered flat connections in the exterior algebra and we note \cite{KM} among other subsequent works. \end{remark}

\end{document}